\documentclass[11pt]{amsart}

\newtheorem{tw}{Theorem}
\newtheorem{lemma}{Lemma}

\newtheorem{prop}{Proposition}

\newcommand{\cal}{\mathcal}

\title{Twistorial construction of generalized K\"ahler manifolds}
\author{Johann Davidov, Oleg Mushkarov}

\address{Institute of Mathematics and Informatics \\
Bulgarian Academy of Sciences\\ Acad. G.Bonchev Str. Bl.8\\
1113 Sofia\\ Bulgaria} \email{jtd@math.bas.bg, muskarov@math.bas.bg}

\begin{document}

\begin{abstract}
The twistor method is applied for obtaining examples of generalized K\"ahler structures
which are not yielded by K\"ahler structures.

 \vspace{0,1cm} \noindent 2000 {\it Mathematics Subject Classification} 53C15, 53C28.

\vspace{0,1cm} \noindent {\it Key words: generalized K\"ahler structures, twistor
spaces}

\end{abstract}

\maketitle \vspace{0.5cm}

\section{Introduction}

The theory of generalized complex structures has been initiated by N. Hitchin
\cite{Hit02} and further developed by M. Gualtieri \cite{Gu}. These structures contain
the complex and symplectic structures as special cases and can be considered as a
complex analog of the notion of a Dirac structure introduced by T. Courant and A.
Weinstein \cite {Cou, CouWei} to unify the Poisson and presymplectic geometries. This
and the fact that the target spaces of supersymmetric $\sigma$-models are generalized
complex manifolds motivate the increasing interest to the generalized complex geometry.

The idea of this geometry is to replace the tangent bundle $TM$ of a smooth manifold
$M$ with the bundle $TM\oplus T^{\ast}M$ endowed with the indefinite metric
$<X+\xi,Y+\eta>=\frac{1}{2}(\xi(Y)+\eta(X))$, $X,Y\in TM$, $\xi,\eta\in T^{\ast}M$. A
generalized K\" ahler structure is, by definition, a pair $\{J_1,J_2\}$ of commuting
generalized complex structures such that the quadratic form $<J_1A,J_2A>$ is positive
definite on $TM\oplus T^{\ast}M$. According to a  result of M. Gualtieri \cite{Gu} the
generalized K\" ahler structures have an equivalent interpretation in terms of the
so-called bi--Hermitian structures.

Any K\" ahler structure yields a generalized K\" ahler structure in a natural way.
Non-trivial examples of such structures can be found in
\cite{ApGaGr,AG,BCG,Hit05,Hit06,K,LT}. The purpose of the present paper is to provide
non-trivial examples of generalized K\" ahler manifolds by means of the R. Penrose
\cite{Pen} twistor construction as developed by M. Atiyah, N. Hitchin and I. Singer
\cite{AtHiSi} in the framework of the Riemannian geometry.

Let $M$ be a $2$-dimensional smooth manifold. Following the general scheme of the
twistor construction we consider the bundle ${\cal P}$ over $M$ whose fibre at a point
$p\in M$ consists of all pairs of commuting generalized complex structures $\{I,J\}$ on
the vector space $T_pM$ such that the form $<IA,JA>$ is positive definite on
$T_pM\oplus T_p^{\ast}M$. The general fibre of ${\cal P}$ admits two natural K\" ahler
structures (in the usual sense) and can be identified in a natural way with the
disjoint union of two copies of the unit bi-disk. Under this identification, the two
structures are defined on the unit bi-disk as $(h\times h,{\cal K}\times (\pm {\cal
K}))$ where $h$ is the Poincare metric on the unit disk and ${\cal K}$ is its standard
complex structure. These two K\"ahler structures yield a generalized K\" ahler
structure on the fibre of ${\cal P}$ according to the Gualtieri result mentioned above.
Moreover, any linear connection $\nabla$ on $M$ gives rise to a splitting of the
tangent bundle $T{\cal P}$ into horizontal and vertical parts and this allows one to
define two commuting generalized almost complex structures ${\cal I}^{\nabla}$ and
${\cal J}^{\nabla}$ on ${\cal P}$ such that the form $<{\cal I}^{\nabla}\cdot,{\cal
J}^{\nabla}\cdot>$ is positive definite on $T{\cal P}\oplus T^{\ast}{\cal P}$. The main
result of the paper states that if the connection $\nabla$ is torsion--free, the
structures ${\cal I}^{\nabla}$ and ${\cal J}^{\nabla}$ are both integrable if and only
if $\nabla$ is flat. Thus any affine structure on $M$ yields a generalized K\" ahler
structure on the $6$-dimensional manifold ${\cal P}$. Note that the only complete
affine $2$-dimensional manifolds are the plane, a cylinder, a Klein bottle, a torus, or
a Mobius band \cite{G,FG}.

\section{Generalized K\"ahler structures}

 Let $W$ be a $n$-dimensional real vector space and $g$ a metric of signature $(p,q)$
on it, $p+q=n$. We shall say that a basis $\{e_1,...,e_{n}\}$ of $W$ is {\it
orthonormal} if $||e_1||^2=...=||e_p||^2=1$, $||e_{p+1}||^2=...=||e_{p+q}||^2=-1$. If
$n=2m$ is an even number and $p=q=m$, the metric $g$ is usually called {\it neutral}.
Recall that a complex structure $J$ on $W$ is called {\it compatible} with the metric
$g$, if the endomorphism $J$ is $g$-skew-symmetric.

\smallskip

   Suppose that $dim\, W=2m$ and  $g$ is of signature $(2p,2q)$, $p+q=m$. Denote by $J(W)$
the set of all complex structures on $W$ compatible with the metric $g$. The group
$O(g)$ of orthogonal transformations of $W$ acts transitively on $J(W)$ by conjugation
and $J(W)$ can be identified with the homogeneous space $O(2p,2q)/U(p,q)$. In
particular, ${\it dim}\,J(W)=m^2-m$. The group $O(2p,2q)$ has four connected
components, while $U(p,q)$ is connected, therefore $J(W)$ has four components.

\smallskip

\noindent {\bf Example 1} (\cite{DM}). The space $O(2,2)/U(1,1)$ is the disjoint union
of two copies of the hyperboloid $x_1^2-x_2^2-x_3^2=1$.

\smallskip

 Consider $J(W)$ as a (closed) submanifold of the vector space $so(g)$ of $g$-skew-symmetric
endomorphisms of $W$. Then the tangent space of $J(W)$ at a point $J$ consists of all
endomorphisms $Q\in so(g)$ anti-commuting with $J$. Thus we have a natural $O(g)$ -
invariant almost complex structure ${\cal K}$ on $J(W)$ defined by ${\cal K}Q=J\circ
Q$. It is easy to check that this structure is integrable.

 Fix an orientation on $W$ and denote by $J^{\pm}(W)$ the set of compatible complex
structures on $W$ that induce $\pm$ the orientation of $W$. The set $J^{\pm}(W)$ has
the homogeneous representation $SO(2p,2q)/U(p,q)$ and, thus, is the union of two
components of $J(W)$.

\smallskip

Suppose that $dim\, W=4$ and $g$ is of split signature $(2,2)$.  Let
$g(a,b)=-\frac{1}{2} Trace (a\circ b)$ be the standard metric of $so(g)$. The
restriction of this metric to the tangent space $T_J$ of $J(W)$ is negative definite
and we set $h=-g$ on $T_J$. Then the complex structure ${\cal K}$ is compatible with
the metric $h$ and $({\cal K},h)$ is a K\"ahler structure on $J(W)$. The space
$J^{\pm}(W)$ can be identified with the hyperboloid  $x_1^2-x_2^2-x_3^2=1$ in ${\Bbb
R}^3$ (see e.g. \cite[Example 5]{DM}) and it is easy to check that, under this
identification, the structure $({\cal K},h)$ on $J^{\pm}(W)$ goes to the standard
K\"ahler structure of the hyperboloid. Thus the Hermitian manifold $(J^{\pm}(W),{\cal
K},h)$ is biholomorphically isometric to the disjoint union of two copies of the unit
disk endowed with the Poincare-Bergman metric (of curvature $-1$).

  Let ${\flat}:T_J\to T_J^{\ast}$ and ${\sharp}={\flat}^{-1}$ be the "musical"
isomorphisms determined by the metric $h$. Denote by $T_J^{\perp}$ the orthogonal
complement of $T_J$ in $so(g)$ with respect to the metric $g$; the space $T_J^{\perp}$
consists of the skew-symmetric endomorphisms of $W$ commuting with $J$. Consider
$T_J^{\ast}$ as the space of linear forms on $so(g)$ vanishing on $T_J^{\perp}$. Then
for every $U\in T_J$ and $\omega\in T_J^{\ast}$ we have $U^{\flat}(A)=-g(U,A)$ and
$g(\omega^{\sharp},A)=-\omega(A)$ for every $A\in so(g)$.

\smallskip

  Now let $V$ be a real vector space and $V^{\ast}$ its dual space. Then
the vector space $V\oplus V^{\ast}$ admits a natural neutral metric defined by
\begin{equation}\label{eq 0.0}
<X+\xi,Y+\eta>=\displaystyle{\frac{1}{2}}(\xi(Y)+\eta(X))
\end{equation}

\smallskip

A {\it generalized complex structure} on the vector space $V$ is, by definition, a
complex structure on the space $V\oplus V^{\ast}$ compatible with its natural neutral
metric \cite{Hit02}. If a vector space $V$ admits a generalized complex structure, it
is necessarily of even dimension \cite{Gu}. We refer to \cite{Gu} for more facts about
the generalized complex structures.

\smallskip

\noindent{\bf Example 2} (\cite{Gu,Hit02,Hit05}). Every complex structure $K$ and every
symplectic form $\omega$ on $V$ (i.e. a non-degenerate $2$-form) induce generalized
complex structures on $V$ in a natural way. If we denote these structures by $J$ and
$S$, respectively, the structure $J$ is defined by $J=K$ on $V$ and $J=-K^{\ast}$ on
$V^{\ast}$, where $(K^{\ast}\xi)(X)=\xi(KX)$ for $\xi\in V^{\ast}$ and $X\in V$.

The map $X\to \imath_X\omega$ (the interior product) is an isomorphism of $V$ onto
$V^{\ast}$. Denote this isomorphism also by $\omega$. Then the structure $S$ is defined
by $S=\omega$ on $V$ and $S=-\omega^{-1}$ on $V^{\ast}$.

\smallskip

\noindent{\bf Example 3} (\cite{Gu,Hit02,Hit05}). Any $2$-form $B\in
\Lambda^{2}V^{\ast}$ acts on $V\oplus V^{\ast}$ via the inclusion
$\Lambda^{2}V^{\ast}\subset \Lambda^{2}(V\oplus V^{\ast})\cong so(V\oplus V^{\ast})$;
in fact this is the action $X+\xi\to \imath_{X}B$;~ $X\in V$, $\xi\in V^{\ast}$. Denote
the latter map again by $B$. Then the invertible map $e^{B}$ is given by $X+\xi\to
X+\xi+\imath_{X}B$ and is an orthogonal transformation of $V\oplus V^{\ast}$. Thus,
given a generalized complex structure $J$ on $V$, the map $e^{B}Je^{-B}$ is also a
generalized complex structure on $V$, called the $B$-transform of $J$.

   Similarly, any $2$-vector $\beta\in \Lambda^{2}V$ acts on $V\oplus V^{\ast}$. If we
identify $V$ with $(V^{\ast})^{\ast}$, so $\Lambda^{2}V\cong
\Lambda^{2}(V^{\ast})^{\ast}$, the action is given by $X+\xi\to \imath_{\xi}\beta\in
V$. Denote this map by $\beta$. Then the exponential map $e^{\beta}$ acts on $V\oplus
V^{\ast}$ via $X+\xi\to X+\imath_{\xi}\beta+\xi$, in particular $e^{\beta}$ is an
orthogonal transformation. Hence, if $J$ is a generalized complex structure on $V$, so
is $e^{\beta}Je^{-\beta}$. It is called the $\beta$-transform of $J$.

\smallskip

   Let $\{e_i\}$ be an arbitrary basis  of $V$ and $\{\eta_i\}$ its
dual basis, $i=1,...,2n$. Then the orientation of the space $V\oplus V^{\ast}$
determined by the basis $\{e_i,\eta_i\}$ does not depend on the choice of the basis
$\{e_i\}$. Further on, we shall always consider $V\oplus V^{\ast}$ with this {\it
canonical orientation}. The sets $J^{\pm}(V\oplus V^{\ast})$ of generalized complex
structures on $V$ inducing $\pm$ the canonical orientation of $V\oplus V^{\ast}$ will
be denoted by $G^{\pm}(V)$.

\smallskip

\noindent {\bf Example 4}. A generalized complex structure on $V$ induced by a complex
structure (see Example 2) always yields the canonical orientation of $V\oplus
V^{\ast}$. A generalized complex structure on $V$ induced by a symplectic form yields
the canonical orientation of $V\oplus V^{\ast}$ if and only if $n=\frac{1}{2}dim\,V$ is
an even number. The $B$- or $\beta$-transform of a generalized complex structure $J$ on
$V$ yields the canonical orientation of $V\oplus V^{\ast}$ if and only if $J$ does so.

\smallskip

\noindent {\bf Example 5}. Let $V$ be a $2$-dimensional real vector space. Take a basis
$\{e_1,e_2\}$ of $V$ and let $\{\eta_1,\eta_2\}$ be its dual basis. Then
$\{Q_1=e_1+\eta_1, Q_2=e_2+\eta_2, Q_3=e_1-\eta_1, Q_4=e_2-\eta_2\}$ is an orthonormal
basis of $V\oplus V^{\ast}$ with respect to the natural neutral metric (\ref{eq 0.0})
and is positively oriented with respect to the canonical orientation of $V\oplus
V^{\ast}$.  Put $\varepsilon_k=||Q_k||^2$, $k=1,...,4$, and define skew-symmetric
endomorphisms of $V\oplus V^{\ast}$ setting $S_{ij}Q_k=\varepsilon_k(\delta_{ik}Q_j -
\delta_{kj}Q_i)$, $1\leq i,j,k\leq 4$. Then the endomorphisms
$$
\begin{array}{lll}
I_1=S_{12}-S_{34}, & & J_1=S_{12}+S_{34},  \\
I_2=S_{13}-S_{24}, & & J_2=S_{13}+S_{24},  \\
I_3=S_{14}+S_{23}, & & J_3=S_{14}-S_{23}
\end{array}
$$
constitute a basis of the space of skew-symmetric endomorphisms of $V\oplus V^{\ast}$.
Let $I\in G^{+}(V)$ and $J\in G^{-}(V)$. Then $I=\sum_{r}x_rI_r$ with
$x_1^2-x_2^2-x_3^2=1$ and $J=\sum_s y_sJ_s$ with $y_1^2-y_2^2-y_3^2=1$. It follows that
$$
\begin{array}{lll}
Ie_1=x_2e_1+(x_1+x_3)e_2, & &  Je_1=y_2e_1+(y_1-y_3)\eta_2,  \\
Ie_2=-(x_1-x_3)e_1-x_2e_2, & & Je_2=y_2e_2-(y_1-y_3)\eta_1,  \\
I\eta_1=-x_2\eta_1+(x_1-x_3)\eta_2, & & J\eta_1=(y_1+y_3)e_2-y_2\eta_1,\\
I\eta_2=-(x_1+x_3)\eta_1+x_2\eta_2, & & J\eta_2=-(y_1+y_3)e_1-y_2\eta_2.
\end{array}
$$
This shows that the restriction of $I$ to $V$ is a complex structure on $V$ inducing
the generalized complex structure $I$ (as in Example 2). In contrast, the generalized
complex structure $J$ is not induced by a complex structure or a symplectic form on
$V$. Moreover $J$ is not a $B$- or $\beta$-transform of such  structures.

{\it A generalized almost complex structure} on an even-dimensional smooth manifold $M$
is, by definition, an endomorphism $J$ of the bundle $TM\oplus T^{\ast}M$ with
$J^2=-Id$ which preserves the natural neutral metric of $TM\oplus T^{\ast}M$. Such a
structure is said to be {\it integrable} or {\it a generalized complex structure} if
its $+i$-eigensubbunle of $(TM\oplus T^{\ast}M)\otimes {\Bbb C}$ is closed under the
Courant bracket \cite{Hit02}. Recall that if $X,Y$ are vector fields on $M$ and
$\xi,\eta$ are $1$-forms, the Courant bracket \cite{Cou} is defined by the formula
$$
[X+\xi,Y+\eta]=[X,Y]+{\cal L}_{X}\eta-{\cal
L}_{Y}\xi-\frac{1}{2}d(\imath_X\eta-\imath_Y\xi),
$$
where $[X,Y]$ on the right hand-side is the Lie bracket and ${\cal L}$ means the Lie
derivative. As in the case of almost complex structures, the integrability condition
for a generalized almost complex structure $J$ is equivalent to the vanishing of its
Nijenhuis tensor $N$, the latter being defined by means of the Courant bracket:
$$
N(A,B)=-[A,B]-J[A,JB]-J[JA,B]+[JA,JB], \> A,B\in TM\oplus T^{\ast}M.
$$

\smallskip

\noindent{\bf Example 6} (\cite{Gu}). A generalized complex structure $K$ induced by an
almost complex structure $K$ on $M$ (see Example 2) is integrable if and only the
structure $K$ is integrable. A generalized complex structure yielded by a
non-degenerate $2$-form $\omega$ on $M$ is integrable if and only if the form $\omega$
is closed.

\smallskip

\noindent{\bf Example 7} (\cite{Gu}).   Let $J$ be a generalized almost complex
structure and $B$ a closed $2$-form on $M$. Then the $B$-transform of $J$,
$e^{B}Je^{-B}$, (see Example 3) is integrable if and only if the structure $J$ is
integrable.

Let us note that the notion of $B$-transform plays an important role in the local
description of the generalized complex structures given by M. Gualtieri \cite{Gu} and
M. Abouzaid - M. Boyarchenko \cite{AbBo}.

\smallskip

The existence of a generalized almost complex structure on a $2n$- dimensional manifold
$M$ is equivalent to the existence of a reduction of the structure group of the bundle
$TM\oplus T^{\ast}M$ to the group $U(n,n)$. Further, to reduce the structure group to
the subgroup $U(n)\times U(n)$ of $U(n,n)$ is equivalent to choosing two commuting
generalized almost complex structures $\{J_1,J_2\}$ such that the quadratic form
$<J_1A,J_2A>$ on $TM\oplus T^{\ast}M$ is positive definite \cite{Gu}. A pair
$\{J_1,J_2\}$ of generalized complex structures with these properties is called an {\it
almost generalized K\"ahler structure}. It is said to be a {\it generalized K\"ahler
structure} if $J_1$ and $J_2$ are both integrable \cite{Gu}.

\smallskip

\noindent{\bf Example 8} (\cite{Gu}). Let $(J,g)$ be a K\"ahler structure on a manifold
$M$ and $\omega$ its K\"ahler form, $\omega(X,Y)=g(JX,Y)$. Let $J_1$ and $J_2$ be the
generalized complex structures on $M$ induced by $J$ and $\omega$. Then the pair
$\{J_1,J_2\}$ is a generalized K\"ahler structure.

\smallskip

\noindent{\bf Example 9} (\cite{Gu}). If $\{J_1,J_2\}$ is a generalized K\"ahler
structure and $B$ is a closed $2$-form, then its $B$-transform
$\{e^{B}J_1e^{-B},e^{B}J_2e^{-B}\}$ is also a generalized K\"ahler structure.

\smallskip

  It has been observed by Gualtieri \cite{Gu} that an almost generalized K\"ahler
structure $\{J_1,J_2\}$ on a manifold $M$ determines the following data on $M$: 1) a
Riemannian metric $g$;~2) two almost complex structures $J_{\pm}$ compatible with
$g$;~3) a $2$-form $b$. Conversely, the almost generalized K\"ahler structure
$\{J_1,J_2\}$ can be reconstructed from the data $(g,J_{+},J_{-},b)$. In fact,
Gualtieri \cite{Gu} has given an explicit formula for $J_1$ and $J_2$ in terms of this
data.

\smallskip

\noindent {\bf Example 10}. Let $V$ be a $2$-dimensional real vector spaces and
$G^{\pm}(V)$ the space of generalized complex structures on $V$ yielding $\pm$ the
canonical orientation of $V\oplus V^{\ast}$. Let $(h,{\cal K})$ be the K\"ahler
structure on $G^{\pm}(V)$ defined above. Consider the manifold $G^{+}(V)\times
G^{-}(V)$ with the product metric $g=h\times h$ and the complex structures $J_{+}={\cal
K}\times {\cal K}$ and $J_{-}={\cal K}\times(-{\cal K})$. According to \cite[formula
(6.3)]{Gu} the generalized K\"ahler structure $\{{\cal I},{\cal J}\}$ on
$G^{+}(V)\times G^{-}(V)$ determined by $g$, $J_{+}$, $J_{-}$ and $b=0$ is given by
\begin{equation}\label{eq GKS-1}
\begin{array}{lll}
{\cal I}(U,V)=I\circ U-V^{\flat}\circ J, & &
{\cal J}(U,V)=J\circ V-U^{\flat}\circ I \\[4pt]
{\cal I}(\varphi,\psi)=-\varphi\circ I+J\circ\psi^{\sharp}, & &{\cal
J}(\varphi,\psi)=-\psi\circ J+I\circ\varphi^{\sharp}
\end{array}
\end{equation}
for $U\in T_{I}G^{+}(V)$, $V\in T_{J}G^{-}(V)$ and $\varphi\in T_{I}^{\ast}G^{+}(V)$,
$\psi\in T_{J}^{\ast}G^{-}(V)$.

\smallskip

  Gualtieri \cite{Gu} has also proved that the integrability condition for
$\{J_1,J_2\}$ can be expressed in terms of the data $(g,J_{+},J_{-},b)$ in a nice way.
In particular, in the case when $b=0$,  the structures $\{J_1,J_2\}$ are integrable if
and only if the almost-Hermitian structures $(g,J_{\pm})$ are Ka\"hlerian.

\smallskip

\noindent {\bf Example 11}. According to the Gualtieri's result the structure $\{{\cal
I},{\cal J}\}$ defined by (\ref{eq GKS-1}) is a generalized K\"ahler structure. Of
course, the integrability of ${\cal I}$ and ${\cal J}$ can be directly proved.

\smallskip

Let $V$ be an even-dimensional real vector space.   The group $GL(V)$ acts on $V\oplus
V^{\ast}$ by letting $GL(V)$ act on $V^{\ast}$ in the standard way. This action
preserves the neutral metric (\ref{eq 0.0}) and the canonical orientation of $V\oplus
V^{\ast}$. Thus, we have an embedding of $GL(V)$ into the group $SO(<~,~>)$ and, via
this embedding, $GL(V)$ acts on the manifold $G^{\pm}(V)$ in a natural manner. Denote
by $P(V)$ the open subset of $G^{+}(V)\times G^{-}(V)$ consisting of those $(I,J)$ for
which the quadratic form $<IA,JA>$ is positive definite on $V\oplus V^{\ast}$. It is
clear that the natural action of $GL(V)$ on $G^{+}(V)\times G^{-}(V)$ leaves $P(V)$
invariant. Suppose that  $dim\, V=2$. Let $I\in G^{+}(V)$ and $J\in G^{-}(V)$. Then it
is easy to see that, under the notations in Example 5, the quadratic form $<IA,JA>$ is
positive definite if and only if either $x_1+x_3>0$, $y_1+y_3>0$ or $x_1+x_3<0$,
$y_1+y_3<0$. This is equivalent to the condition that either $x_1>0$, $y_1>0$ or
$x_1<0$, $y_1<0$. Thus $P(V)$ is the disjoint union of two products of one-sheeted
hyperboloids. Therefore $P(V)$ endowed with the complex structure ${\cal K}\times {\cal
K}$ and the metric $h\times h$ is biholomorphically isometric to the disjoint union of
two copies of the unit bi-disk endowed with the Bergman metric. Note also that, when
$dim\, V=2$, every $I\in G^{+}(V)$ {\it commutes} with every $J\in G^{-}(V)$ (see
Example 5). Thus, in this case, every pair $(I,J)\in P(V)$ is a generalized K\"ahler
structures on the manifold $V$.

\section{The twistor space of generalized K\"ahler structures}

 Let $M$ be a smooth manifold of dimension $2$. Denote
by $\pi:{\cal G}^{\pm}\to M$ the bundle over $M$ whose fibre at a point $p\in M$
consists of all generalized complex structures on $T_pM$ that induce $\pm$ the
canonical orientation of $T_pM\oplus T^{\ast}_pM$. This is the associated bundle
$$
GL(M)\times_{GL(2,{\Bbb R} )} G^{\pm}({\Bbb R}^{2}),
$$
where $GL(M)$ denotes the principal bundle of linear frames on $M$. Consider the
product bundle $\pi:{\cal G}^{+}\times {\cal G}^{-}\to M$ and denote by ${\cal P}$ its
open subset consisting of those pairs $K=(I,J)$ for which the quadratic form $<IA,JA>$
on $T_pM\oplus T_p^{\ast}M$, $p=\pi(K)$, is positive definite. Clearly ${\cal P}$ is
the associated bundle
$$
{\cal P}=GL(M)\times_{GL(2,{\Bbb R} )} P({\Bbb R}^{2}).
$$
The projection maps of the bundles  ${\cal G}^{\pm}$ and ${\cal P}$ to the base space
$M$ will be denoted by $\pi$.

 Let $\nabla$ be a linear connection on $M$. Following the standard twistor
construction we can define two commuting almost generalized complex structures ${\cal
I}^{\nabla}$ and ${\cal J}^{\nabla}$ on ${\cal P}$ as follows: The connection $\nabla$
gives rise to a splitting ${\cal V}\oplus {\cal H}$ of the tangent bundle of any bundle
associated to $GL(M)$ into vertical and horizontal parts. The vertical space ${\cal
V}_K$ of ${\cal P}$ at a point $K=(I,J)$ is the direct sum ${\cal V}_K={\cal V}_I{\cal
G}^{+}\oplus {\cal V}_J{\cal G}^{-}$ of vertical spaces and we define  ${\cal
I}^{\nabla}$ and ${\cal J}^{\nabla}$ on ${\cal V}_K$ by means of (\ref{eq GKS-1}) where
the "musical" isomorphisms are determined by the metric $h$ on ${\cal V}_J{\cal G}^{+}$
and ${\cal V}_J{\cal G}^{-}$.

The horizontal space ${\cal H}_K$ is isomorphic via the differential $\pi_{\ast K}$ to
the tangent space $T_{p}M, p=\pi(K)$. Denoting $\pi_{\ast K}|{\cal H}$ by $\pi_{{\cal
H}}$, we define ${\cal I}^{\nabla}$ and ${\cal J}^{\nabla}$ on ${\cal H}_K\oplus {\cal
H}_K^{\ast}$ as the lift of the endomorphisms $I$ and $J$ by the map $\pi_{{\cal
H}}\oplus (\pi_{{\cal H}}^{-1})^{\ast}$.

\smallskip

\noindent

{\bf Remark}.  Neither of the generalized almost complex structures ${\cal I}^{\nabla}$
and ${\cal J}^{\nabla}$ is induced by an almost complex or symplectic structure on
${\cal P}$. Moreover they are not $B$- or $\beta$-transforms of such structures.

\smallskip

     Further on, the generalized almost complex structures ${\cal I}^{\nabla}$ and
${\cal J}^{\nabla}$ will be simply denoted by ${\cal I}$ and ${\cal J}$ when the
connection $\nabla$ is understood. The image of every $A\in T_pM\oplus T_p^{\ast}M$
under the map $\pi_{{\cal H}}^{-1}\oplus \pi_{{\cal H}}^{\ast}$ will be denoted by
$A^h$. The elements of ${\cal H}_J^{\ast}$, resp. ${\cal V}_J^{\ast}$, will be
considered as $1$-forms on $T_J{\cal G}$ vanishing on ${\cal V}_J$, resp. ${\cal H}_J$.

\smallskip

Let $K=(I,J)\in{\cal P}$, $A\in T_{\pi(K)}M\oplus T_{\pi(K)}^{\ast}M$, $W=(U,V)\in{\cal
V}_K$ and $\Theta=(\varphi,\psi)\in{\cal V}_{K}^{\ast}$. Then we have
$$
<{\cal I}(A^h+W+\Theta),{\cal
J}(A^h+W+\Theta)>=<IA,JA>+||U||_h^2+||V||_h^2+||\varphi||_h^2+||\psi||_h^2.
$$
Therefore the quadratic form $<{\cal I}\cdot,{\cal J}\cdot>$ is positive definite. Thus
the pair $({\cal I},{\cal J})$ is an almost generalized K\"ahler structure.

\smallskip

We shall show that for a torsion-free connection $\nabla$ the integrability condition
for ${\cal I}$ and ${\cal J}$ can be expressed in terms of the curvature of $\nabla$
(as is usual in the twistor theory).

\smallskip

Let $A(M)$ be the bundle of the endomorphisms of $TM\oplus T^{\ast}M$ which are skew-
symmetric with respect to its natural neutral metric $<~,~>$; the fibre of this bundle
at a point $p\in M$ will be denoted by $A_p(M)$. The connection $\nabla$ on $TM$
induces a connection on $A(M)$, thus a connection on the bundle $A(M)\oplus A(M)$, both
denoted again by $\nabla$.

  Consider the bundle ${\cal P}$ as a subbundle of the bundle
$\pi: A(M)\oplus A(M)\to M$. Then the inclusion of ${\cal P}$ is fibre-preserving and
the horizontal space of ${\cal P}$ at a point $K$ coincides with the horizontal space
of $A(M)\oplus A(M)$ at that point since the inclusion $P({\Bbb R}^{2})\subset
so(2,2)\times so(2,2)$ is $SO(2,2)$-equivariant.

Let $(U,x_1, x_{2})$ be a local coordinate system of $M$ and $\{Q_1,...,Q_{4}\}$ an
orthonormal frame of $TM\oplus T^{\ast}M$ on $U$. Set $\varepsilon_k=||Q_k||^2$,
$k=1,...,4$, and define sections $S_{ij}$, $1\leq i,j\leq {4}$, of $A(M)$ by the
formula
\begin{equation}\label{eq Sij}
S_{ij}Q_k=\varepsilon_k(\delta_{ik}Q_j - \delta_{kj}Q_i).
\end{equation}
Then $S_{ij}$, $i<j$, form an orthogonal frame of $A(M)$ with respect to the metric
$<a,b>=\displaystyle{-\frac{1}{2}}Trace\,(a\circ b);\, a,b\in A(M)$; moreover
$||S_{ij}||^2=\varepsilon_i\varepsilon_j$ for $i\neq j$.  For $c=(a,b)\in A(M)\oplus
A(M)$, we set
$$\tilde x_{m}(c)=x_{m}\circ\pi(c),~ y_{ij}(c)=\varepsilon_i\varepsilon_j<a,S_{ij}>,
~ z_{ij}(c)=\varepsilon_i\varepsilon_j<b,S_{ij}>.$$ Then $(\tilde x_{m},y_{ij},
z_{kl})$, $m=1,2$, $1\leq i < j\leq 4$, $1\leq k < l\leq 4$, is a local coordinate
system on the total space of the bundle $A(M)\oplus A(M)$. Note that $(\tilde
x_{m},y_{ij})$ and $(\tilde x_{m},z_{kl})$ are local coordinate systems of the manifold
$A(M)$.

Let
$$
U=\sum_{i<j}u_{ij}\frac{\partial}{\partial y_{ij}}(I),\quad
V=\sum_{i<j}v_{ij}\frac{\partial}{\partial z_{ij}}(J)
$$
be vertical vectors of ${\cal G}^{+}$ and ${\cal G}^{-}$ at some points $I$ and $J$
with $\pi(I)=\pi(J)$. It is convenient to set $u_{ij}=-u_{ji}$, $v_{ij}=-v_{ji}$ for
$i\geq j$, $1\leq i,j\leq {4}$. Then the endomorphism $U$ of $T_{p}M\oplus
T_{p}^{\ast}M$, $p=\pi(I)$, is determined by $UQ_i=\sum_{j=1}^{4}
\varepsilon_iu_{ij}Q_j$; similarly for the endomorphism $V$ of $T_{p}M\oplus
T_{p}^{\ast}M$.  Moreover
$$
{\cal K}_{I}^{\ast}U^{\flat}=-(IU)^{\flat}=
\sum_{i<j}\varepsilon_i\varepsilon_j\sum_{r=1}^{4}u_{ir}y_{rj}(I)\varepsilon_r(dy_{ij})_I.
$$
Similar formula holds for ${\cal K}_{J}^{\ast}V^{\flat}$. Thus we have
\begin{equation}\label{cal I/ver}
{\cal I}(U,V)=\sum_{i<j}\sum_{r}u_{ir}y_{rj}(I)\varepsilon_{r}\frac{\partial}{\partial
y_{ij}}(I) -
\sum_{k<l}\varepsilon_{k}\varepsilon_{l}\sum_{s}v_{ks}z_{sl}(J)\varepsilon_{s}(dz_{kl})_J
\end{equation}
and
\begin{equation}\label{cal J/ver}
{\cal J}(U,V)=\sum_{k<l}\sum_{s}v_{ks}z_{sl}(J)\varepsilon_{s}\frac{\partial}{\partial
z_{kl}}(J) -
\sum_{i<j}\varepsilon_{i}\varepsilon_{j}\sum_{r}u_{ir}y_{rj}(I)\varepsilon_{r}(dy_{ij})_I.
\end{equation}

 Note also that, for every $A\in T_{p}M\oplus T_{p}^{\ast}M$, we have
\begin{equation}\label{A^h}
A^h=\sum_{i=1}^{4n}(<A,Q_i>\circ\pi)\varepsilon_i Q_i^h
\end{equation}
and
\begin{equation}\label{cal I/hor}
{\cal I}A^h=\sum_{i,j=1}^{4}(<A,Q_i>\circ\pi)y_{ij}Q_j^h, \quad {\cal
J}A^h=\sum_{k,l=1}^{4}(<A,Q_k>\circ\pi)z_{kl}Q_l^h.
\end{equation}

   For each vector field
$$X=\sum_{i=1}^{2} X^{i}\frac{\partial}{\partial x_i}$$
on $U$, the horizontal lift $X^h$ on $\pi^{-1}(U)$ is given by
\begin{equation}\label{eq 3.1}
\begin{array}{lll}
X^{h}&=&\displaystyle{\sum_{m}(X^{m}\circ\pi)\frac{\partial}{\partial\tilde x_m}}\\
                                                    \\
& &-\displaystyle{\sum_{i<j}\sum_{a<b}
y_{ab}(<\nabla_{X}S_{ab},S_{ij}>\circ\pi)\varepsilon_i\varepsilon_j\frac{\partial}{\partial
y_{ij}}} \\
                                                        \\
& &- \displaystyle{\sum_{k<l}\sum_{c<d}
z_{cd}(<\nabla_{X}S_{cd},S_{kl}>\circ\pi)\varepsilon_k\varepsilon_l\frac{\partial}{\partial
z_{kl}}}.
\end{array}
\end{equation}

     Let $c=(a,b)\in A(M)\oplus A(M)$ and $p=\pi(c)$. Then (\ref{eq 3.1}) implies
that, under the standard identification of $T_{c}(A_p(M)\oplus A_p(M))$ with the vector
space $A_p(M)\oplus A_p(M)$, we have
\begin{equation}\label{eq 3.2}
[X^{h},Y^{h}]_{c}=[X,Y]^h_c + R(X,Y)c,
\end{equation}
where $R(X,Y)c=(R(X,Y)a,R(X,Y)b)$ is the curvature of the connection $\nabla$ on
$A(M)\oplus A(M)$ (for the curvature tensor we adopt the following definition:
$R(X,Y)=\nabla_{[X,Y]}- [\nabla_{X},\nabla_{Y}]$).

\smallskip

 {\it Notation}. Let $K=(I,J)\in {\cal P}$ and $p=\pi(K)$. There exists an oriented
orthonormal basis $\{a_1,...,a_{4}\}$ of $T_pM\oplus T_{p}^{\ast}M$ such that
$a_2=Ia_1$, $a_4=Ia_3$ and $Ja_1=\varepsilon a_2$, $ Ja_3=-\varepsilon a_4$, where
$\varepsilon=+1$ or $-1$.  Let $\{Q_i\}$, $i=1,...,4$, be an oriented orthonormal frame
of $TM\oplus T^{\ast}M$ near the point $p$ such that
$$
Q_i(p)=a_i \mbox { and } \nabla Q_i|_p=0, \> i=1,...,4.
$$
Define sections $S$ and $T$ of $A(M)$ by setting
$$
\begin{array}{lll}
SQ_1=Q_2,\quad SQ_2=-Q_1,\quad SQ_3=Q_4,\quad SQ_4=-Q_3 \\
TQ_1=\varepsilon Q_2,\quad JQ_2=-\varepsilon Q_1,\quad TQ_3=-\varepsilon Q_4,\quad
TQ_4=\varepsilon Q_3.
\end{array}
$$
Then $\nu=(S,T)$ is a section of ${\cal P}$ such that
$$
\nu(p)=K, \quad \nabla\nu|_p=0
$$
(considering $\nu$ as a section of $A(M)\oplus A(M)$). Thus $X^h_K=\nu_{\ast}X$ for
every $X\in T_pM$.

 Further, given a smooth manifold $N$, the natural projections of $TN\oplus T^{\ast}N$
onto $TN$ and $T^{\ast}N$ will be denoted by $\pi_1$ and $\pi_2$, respectively.

  We shall use the above notations throughout this section.

\vspace{0.1cm}

The next three technical lemmas can be easily proved by means of  (\ref{cal I/hor}),
(\ref{eq 3.1}) and (\ref{eq 3.2}).

\begin{lemma}\label {brackets}
If $A$ and $B$ are sections of the bundle $TM\oplus T^{\ast}M$ near $p$, then:
\begin{enumerate}
\item[$(i)$]
$[\pi_1(A^h),\pi_1({\cal I }B^h)]_K=[\pi_1(A),\pi_1(SB)]^h_K+R(\pi_1(A),\pi_1(IB))K.$
\vspace{0.2cm}
\item[$(ii)$]
$[\pi_1({\cal I}A^h),\pi_1({\cal
I}B^h)]_K=[\pi_1(SA),\pi_1(SB)]^h_K+R(\pi_1(IA),\pi_1(IB))K.$
\end{enumerate}
\end{lemma}

\smallskip

\begin{lemma}\label {Lie deriv}
Let $A$ and $B$ be sections of the bundle $TM\oplus T^{\ast}M$ near $p$, and let $Z\in
T_pM$, $W=(U,V)\in {\cal V}_K={\cal V}_I{\cal G}^{+}\oplus {\cal V}_J{\cal G}^{-}$.
Then:
\begin{enumerate}
\item[$(i)$]
\hspace{0.2cm}$({\cal L}_{\pi_1(A^h)}{\pi_2(B^h)})_K=({\cal
L}_{\pi_1(A)}{\pi_2(B)})^h_K.$ \vspace{0.2cm}
\item[$(ii)$]
\hspace{0.2cm}$({\cal L}_{\pi_1(A^h)}{\pi_2({\cal I}B^h)})_K=({\cal
L}_{\pi_1(A)}{\pi_2(SB)})^h_K.$ \vspace{0.2cm}
\item[$(iii)$] \ \\
$
\begin{array}{lll}
({\cal L}_{\pi_1({\cal I}A^h)}\pi_2(B^h))_K(Z^h+W)=\\
                                                    \\
({\cal L}_{\pi_1(SA)}\pi_2(B))^h_K(Z^h)+(\pi_2(B))_p(\pi_1(UA)).
\end{array}
$
 \vspace{0.2cm}
\item[$(iv)$] \ \\
$
\begin{array}{lll}
({\cal L}_{\pi_1({\cal I}A^h)}\pi_2({\cal I}B^h))_K(Z^h+W)=\\
                                                             \\
({\cal L}_{\pi_1(SA)}\pi_2(SB))^h_K(Z^h)+(\pi_2(IB))_p(\pi_1(UA)).
\end{array}
$
\end{enumerate}
\end{lemma}

\smallskip

\begin{lemma}\label {Half-diff}
Let $A$ and $B$ are sections of the bundle $TM\oplus T^{\ast}M$ near $p$. Let $Z\in
T_pM$ and $W=(U,V)\in {\cal V}_K={\cal V}_I{\cal G}^{+}\oplus {\cal V}_J{\cal G}^{-}$.
Then:
\begin{enumerate}
\item[$(i)$]
\hspace{0.2cm}$(d\>\imath_{\pi_1(A^h)}\pi_2(B^h))_K=(d\>\imath_{\pi_1(A)}\pi_2(B))^h_K$
\vspace{0.2cm}
\item[$(ii)$] \ \\
$
\begin{array}{lll}
(d\>\imath_{\pi_1(A^h)}\pi_2({\cal I}B^h))_K(Z^h+W)= \\
                                                                \\
(d\>\imath_{\pi_1(A)}\pi_2(SB))^h_K(Z^h)+ (\pi_2(UB))_p(\pi_1(A))
\end{array}
$ \vspace{0.2cm}
\item[$(iii)$] \ \\
$
\begin{array}{lll}
(d\>\imath_{\pi_1({\cal I}A^h)}\pi_2(B^h))_K(Z^h+W)= \\
                                                               \\
(d\>\imath_{\pi_1(SA)}\pi_2(B))^h_K(Z^h)+ (\pi_2(B))_p(\pi_1(UA))
\end{array}
$ \vspace{0.2cm}
\item[$(iv)$] \ \\
$
\begin{array}{lll}
(d\>\imath_{\pi_1({\cal I}A^h)}\pi_2({\cal I}B^h))_K(Z^h+W)= \\
                                                                                \\
(d\>\imath_{\pi_1(SA)}\pi_2(SB))^h_K(Z^h)+ (\pi_2(UB))_p(\pi_1(IA))+
(\pi_2(IB))_p(\pi_1(UA)
\end{array}
$

\end{enumerate}
\end{lemma}

\smallskip

\begin{prop}\label {Nijenhuis-hor}
Suppose that the connection $\nabla$ is torsion-free and let $K=(I,J)\in {\cal P}$.
Then
\begin{enumerate}
\item[$(i)$] \quad
$N^{{\cal I}}(A^h,B^h)=0$ for  every  $A,\> B\in T_{\pi(K)}M\oplus T_{\pi(K)}^{\ast}M$.
          \\
\item[$(ii)$] \quad
$N^{{\cal J}}(A^h,B^h)=0$ for  every  $A,\> B\in T_{\pi(K)}M\oplus T_{\pi(K)}^{\ast}M$
if and only if $R(X,Y)J=0$ for every $X,Y\in T_{\pi(K)}M$.
\end{enumerate}
\end{prop}

\begin{proof}
First we shall show that
\begin{equation}\label{eq Nij-I}
\begin{array}{lll}
N^{{\cal I}}(A^h,B^h)_K =\\[4pt] -R(\pi_1(A),\pi_1(B))I-I\circ
R(\pi_1(A),\pi_1(IB))I\\[4pt]
-I\circ R(\pi_1(IA),\pi_1(B))I+R(\pi_1(IA),\pi_1(IB))I\\[4pt]
 -R(\pi_1(A),\pi_1(B))J+R(\pi_1(IA),\pi_1(IB))J\\[4pt]
 +{\cal K}_{J}^{\ast}(R(\pi_1(A),\pi_1(IB)J)^{\flat}+ {\cal
K}_{J}^{\ast}(R(\pi_1(IA),\pi_1(B)J)^{\flat}.
\end{array}
\end{equation}
Similar formula holds for the Nijenhuis tensor $N^{\cal J}$ with interchanged roles of
$I$ and $J$ in the right-hand side of (\ref{eq Nij-I}).

 Set $p=\pi(K)$ and extend  $A$ and $B$ to (local) sections of $TM\oplus T^{\ast}M$,
denoted again by $A,B$, in such a way that $\nabla A|_p=\nabla B|_p=0$.

 Let $\nu=(S,T)$ be the section of ${\cal P}$ defined above with the property that
$\nu(p)=K$ and $\nabla \nu|_p=0$  ($\nu$ being considered as a section of $A(M)\oplus
A(M)$).

 According to Lemmas~\ref{brackets}, \ref{Lie deriv} and \ref{Half-diff}, the
part of $N^{{\cal I}}(A^h,B^h)_K$ lying in  ${\cal H}_K\oplus {\cal H}^{\ast}_K$ is
given by
\begin{equation}\label{eq hor-part Nij-I}
\begin{array}{c}
({\cal H}\oplus {\cal H}^{\ast})N^{{\cal I}}
(A^h,B^h)_K=\\[4pt]
(-[A,B]-S[A,SB]-S[SA,B]+[SA,SB])^h_K.
\end{array}
\end{equation}

Note that we have $\nabla\pi_1(A)|_p=\pi_1(\nabla A|_p)=0$  and  $\nabla\pi_1(SA)|_p=
\pi_1((\nabla S)|_p(A)\\+ S(\nabla A|_p))=0$. Similarly, $\nabla\pi_2(A)|_p=0$ and
$\nabla\pi_2(SA)|_p=0$. We also have $\nabla\pi_1(B)|_p=0$, $\nabla\pi_1(SB)|_p=0$ and
$\nabla\pi_2(B)|_p=0$, $\nabla\pi_2(SB)|_p=0$. Now, since $\nabla$ is torsion-free, we
can easily see that every bracket in (\ref{eq hor-part Nij-I}) vanishes by means of the
following simple observation: Let $Z$ be a vector field and $\omega$ a $1$-form on $M$
such that $\nabla Z|_p=0$ and $\nabla\omega|_p=0$. Then for every $T\in T_pM$
$$
({\cal L}_{Z}\omega)(T)_p =(\nabla_{Z}\omega)(T)_p=0 ~\mbox { and
}~(d\,\imath_{Z}\omega)(T)_p=(\nabla_{T}\omega)(Z)_p=0.
$$

   By Lemmas~\ref{brackets} -- \ref{Half-diff}, the
part of $N^{{\cal I}}(A^h,B^h)_K$ lying in ${\cal V}_K$ is
$$
\begin{array}{lll}
 -R(\pi_1(A),\pi_1(B))I-I\circ R(\pi_1(A),\pi_1(IB))I\\[4pt]
-I\circ R(\pi_1(IA),\pi_1(B))I+R(\pi_1(IA),\pi_1(IB))I\\[4pt]
 -R(\pi_1(A),\pi_1(B))J+R(\pi_1(IA),\pi_1(IB))J
 \end{array}
$$

Finally, the part of $N^{{\cal I}}(A^h,B^h)_K$ lying in ${\cal V}_K^{\ast}$ is the
vertical form whose value at every vertical vector $W=(U,V)\in {\cal V}_K$ is equal to
$$
\begin{array}{lll}
\frac{1}{2}\{-\pi_2(IUB)(\pi_1(A))-\pi_2(A)(\pi_1(IUB))\\[4pt]
+\pi_2(IUA)(\pi_1(B))+\pi_2(B)(\pi_1(IUA))\\[4pt]
+\pi_2(IB)(\pi_1(UA))+\pi_2(UA)(\pi_1(IB))\\[4pt]
-\pi_2(IA)(\pi_1(UB)) -\pi_2(UB)(\pi_1(IA)) \}\\[4pt]
+{\cal K}_{J}^{\ast}(R(\pi_1(A),\pi_1(IB)J)^{\flat}+ {\cal
K}_{J}^{\ast}(R(\pi_1(IA),\pi_1(B)J)^{\flat}.
\end{array}
$$
The endomorphism $U$ of $T_pM\oplus T_p^{\ast}M$ is skew-symmetric with respect to the
metric $<~,~>$ \, and anti-commutes with $I$. Thus we have
$$
<IUA,B>=<IA,UB>.
$$
This identity reads as
$$
\pi_2(IUA)(\pi_1(B))+\pi_2(B)(\pi_1(IUA))=\pi_2(IA)(\pi_1(UB)) +\pi_2(UB)(\pi_1(IA)).
$$
Therefore the part of $N^{{\cal I}}(A^h,B^h)_K$ lying in ${\cal V}_K^{\ast}$ is
$$
{\cal K}_{J}^{\ast}(R(\pi_1(A),\pi_1(IB)J)^{\flat}+ {\cal
K}_{J}^{\ast}(R(\pi_1(IA),\pi_1(B)J)^{\flat}.
$$

This proves formula (\ref{eq Nij-I}).

 Now let $\{Q_1, Q_2=IQ_1, Q_3, Q_4=IQ_3\}$ be an orthonormal basis
of $T_pM\oplus T_p^{\ast}M$. To prove that $N^{{\cal I}}(A^h,B^h)_K=0$ it is enough to
show that $N^{{\cal I}}(Q_1^h,Q_3^h)_K=0$ since $N^{{\cal I}}({\cal I}E,F)=N^{{\cal I
}}(E,{\cal I}F)=-{\cal I}N^{{\cal I}}(E,F)$ for every $E,F\in T{\cal P}$.

Let $\pi_1(Q_i)=e_i$, $i=1,...,4$. Then, according to (\ref{eq Nij-I})
$$
\begin{array}{lll}
N^{{\cal
I}}(Q_1^h,Q_3^h)&=&[-R(e_1,e_3)I+R(e_2,e_4)I]-I\circ[R(e_1,e_4)I+R(e_2,e_3)I]\\[4 pt]
& & + {\cal K}_J^{\ast}(R(e_1,e_4)J+R(e_2,e_3)J)^{\flat}.
\end{array}
$$
Since $I$ yields the canonical orientation of  $T_pM\oplus T_p^{\ast}M$, the latter
expression vanishes in view of  the following simple algebraic fact proved in
\cite{DM}:

\begin{lemma}\label {Orient}
 Let $V$ be a $2$-dimensional real vector space and let $\{Q_i=e_i+\eta_i\}$, $1\leq i\leq 4$,
be an orthonormal basis of the space $V\oplus V^{\ast}$ endowed with its natural
neutral metric {\rm {(\ref{eq 0.0})}}. Then $\{e_1,e_2\}$ is a bases of $V$ and
$$
\begin{array}{lll}
e_3=a_{11}e_1+a_{12}e_2 \\
e_4=a_{21}e_1+a_{22}e_2
\end{array}
$$
where $A=[a_{kl}]$ is an orthogonal matrix. If $det\,A=1$, the basis $\{Q_i\}$ yields
the canonical orientation of $V\oplus V^{\ast}$ and if $det\,A=-1$ it yields the
opposite one.
\end{lemma}

  To prove statement $(ii)$, take an orthonormal basis $\{\bar Q_1, \bar Q_2=J\bar Q_1,
\bar Q_3, \bar Q_4=J\bar Q_3\}$ and set $\pi_1(\bar Q_i)=e_i$, $i=1,...,4$. Suppose
that $N^{{\cal J}}(\bar Q_1^h,\bar Q_3^h)=0$. Then, according to the analog of (\ref{eq
Nij-I}) for $N^{{\cal J}}(A^h,B^h)_K$, we have
$$
-R(e_1,e_3)J+R(e_2,e_4)J-J\circ[R(e_1,e_4)J+R(e_2,e_3)J]=0
$$
Since $J$ yields the orientation of $T_pM\oplus T_p^{\ast}M$ opposite to the canonical
one, then, by Lemma~\ref{Orient}, $e_3=\cos{t}\,e_1+\sin{t}\,e_2$, $e_4=\sin{t}\,
e_1-\cos{t}\,e_2$ for some $t\in {\Bbb R}$. Thus
$$
-sin{t}\cdot R(e_1,e_2)J +cos{t}\cdot J\circ R(e_1,e_2)J=0,
$$
which implies
$$
cos{t}\cdot R(e_1,e_2)J +sin{t}\cdot J\circ R(e_1,e_2)J=0.
$$
Therefore $R(e_1,e_2)J=0$, so $R(X,Y)J=0$ for every $X,Y\in T_pM$.

Conversely, if the latter identity holds, the analog of (\ref{eq Nij-I}) shows that
$N^{{\cal J}}(A^h,B^h)_K=0$.
\end{proof}

\smallskip

\begin{prop}\label {Nijenhuis-hor-ver}
Suppose that the connection $\nabla$ is torsion-free and let $K=(I,J)\in {\cal P}$,
Then
\begin{enumerate}
\item[$(i)$] \quad
$N^{{\cal I}}(A^h,W)=0$ for  every  $A\in T_{\pi(K)}M\oplus T_{\pi(K)}^{\ast}M$ and
$W\in{\cal V}_K$ if and only if $R(X,Y)J=0$ for every $X,Y\in T_{\pi(K)}M$.
\
          \\
\item[$(ii)$] \quad
$N^{{\cal J}}(A^h,W)=0$ for  every  $A\in T_{\pi(K)}M\oplus T_{\pi(K)}^{\ast}M$ and
$W\in{\cal V}_K$ if and only if $R(X,Y)I=0$ for every $X,Y\in T_{\pi(K)}M$.
\end{enumerate}
\end{prop}

\begin{proof}
Set $p=\pi(K)$ and $W=(U,V)$. Extend $A$ to a section of $TM\oplus T^{\ast}M$  denoted
again by $A$. Take sections $a$ and $b$ of $A(M)$ such that
$$
a(p)=U, \quad b(p)=V, \quad \nabla a|_p=\nabla b|_p=0.
$$
Define vertical vector fields $\widetilde a$ and $\widetilde b$ on ${\cal G}^{+}$ and
${\cal G}^{-}$, respectively, setting
\begin{equation}\label{eq def-ver-ab}
\widetilde a_{I^{'}}=a_{\pi(I^{'})}+I^{'}\circ a_{\pi(I^{'})}\circ I^{'}, ~I^{'}\in
{\cal G}^{+} \mbox { and } \widetilde b_{J^{'}}=b_{\pi(J^{'})}+J^{'}\circ
b_{\pi(J^{'})}\circ J^{'}, ~ J^{'}\in {\cal G}^{-}.
\end{equation}
Then
$$
\widetilde{W}_{(I^{'},J^{'})}=(\widetilde a_{I^{'}},\widetilde b_{J^{'}}), \quad
(I^{'},J^{'})\in {\cal P},
$$
is a vertical vector field on ${\cal P}$ with $\widetilde{W}_K=2W$.

  Let $a(Q_i)=\sum_{j}\varepsilon_i a_{ij}Q_j$,\quad $b(Q_i)=\sum_{j}\varepsilon_i b_{ij}Q_j$.
Then, in the local coordinates introduced above,
\begin{equation}\label{eq W-tilde}
\widetilde W=\sum_{i<j}(\widetilde a_{ij}\frac{\partial}{\partial y_{ij}}+\widetilde
b_{ij}\frac{\partial}{\partial z_{ij}}),
\end{equation}
where
$$
\widetilde
a_{ij}=a_{ij}\circ\pi+\sum_{k,l}y_{ik}(a_{kl}\circ\pi)y_{lj}\varepsilon_k\varepsilon_l,
\quad \widetilde
b_{ij}=b_{ij}\circ\pi+\sum_{k,l}z_{ik}(b_{kl}\circ\pi)z_{lj}\varepsilon_k\varepsilon_l.
$$
In view of (\ref{eq 3.1}), for any vector field $X$ on $M$ near the point $p$, we have
\begin{equation}\label{eq Xh-y-z}
X^h_K=\sum_{m}X^{m}(p)\frac{\partial}{\partial\tilde x_m}(K),\quad
[X^h,\frac{\partial}{\partial y_{ij}}]_K=[X^h,\frac{\partial}{\partial z_{ij}}]_K=0,
\end{equation}
and
$$
0=(\nabla_{X_{p}}a)(Q_i)=\sum_{j}\varepsilon_iX_{p}(a_{ij})Q_j,\quad
0=(\nabla_{X_{p}}b)(Q_i)=\sum_{j}\varepsilon_iX_{p}(b_{ij})Q_j
$$
since $\nabla Q_i|_p=0$ and $\nabla S_{ij}|_p=0$. In particular,
$X_p(a_{ij})=X_p(b_{ij})=0$, hence
\begin{equation}\label{eq der a-b}
X^h_K(\widetilde a_{ij})=X^h_K(\widetilde b_{ij})=0.
\end{equation}
Now simple calculations making use of (\ref{eq Xh-y-z}), (\ref{eq 3.1}) and (\ref{eq
W-tilde}) give
\begin{equation}\label{eq hor-ver-tilde}
[X^h,\widetilde W]_K=0.
\end{equation}
Let $\omega$ be a $1$-form on $M$. It is easy to see that for every vertical vector
field $W^{'}$ on ${\cal P}$
\begin{equation}\label{eq form-ver-tilde}
[\omega^h,W^{'}]=0.
\end{equation}
Therefore, by (\ref{eq hor-ver-tilde}) and (\ref{eq form-ver-tilde}), we have
\begin{equation}\label{eq bra-1}
[A^h,\widetilde W]_K=0.
\end{equation}
Next, in view of (\ref{eq form-ver-tilde}),  (\ref{cal I/ver}), (\ref{eq Xh-y-z}) and
(\ref{eq der a-b}), we have
$$
[A^h,{\cal I}\widetilde W]_K=[\pi_1(A^h),{\cal I}\widetilde W]_K=({\cal
L}_{\pi_1(A^h)}\pi_2({\cal I}\widetilde W))_K.
$$
Let $W^{'}=(U^{'},V^{'})\in {\cal V}_K$. Take sections $a^{'}, b^{'}$ of $A(M)$ such
that $a^{'}(p)=U^{'}$, $b^{'}(p)=V^{'}$, $\nabla a^{'}|_p=\nabla b^{ '}|_p=0$. Define
vertical vector fields $\widetilde {a^{'}}$ and $\widetilde {b^{'}}$ on ${\cal G}^{+}$
and ${\cal G}^{-}$ by means of (\ref{eq def-ver-ab}) and set $\widetilde
{W^{'}}=(\widetilde {a^{'}},\widetilde {b^{'}})$ on ${\cal P}$. Then $[X^h,\widetilde
{W^{'}}]_K=0$ for every vector field $X$ near the point $p$ and an easy computation
making use of (\ref{cal I/ver}), (\ref{eq Xh-y-z}) and (\ref{eq der a-b}) gives
$$
({\cal L}_{\pi_1(A^h)}\pi_2({\cal I}\widetilde W))_K(W^{'})=\frac{1}{2}({\cal
L}_{\pi_1(A^h)}\pi_2({\cal I}\widetilde W))_K(\widetilde {W^{'}})=0.
$$
Moreover, for every vector field $Z$ on $M$ near the point $p$ we have
$$
\begin{array}{c} ({\cal L}_{\pi_1(A^h)}\pi_2({\cal I}\widetilde W))_K(Z^h)=-\pi_2({\cal
I}\widetilde W)([\pi_1(A^h),Z^h]_K)=\\[4pt]
2V^{\flat}(J\circ R(\pi_1(A),Z)J)=2<JV,R(\pi_1(A),Z)J>.
\end{array}
$$
by (\ref{eq GKS-1}) and (\ref{eq 3.2}). It is convenient to define a $1$-form
$\gamma_A$ on $T_pM$ setting
$$\gamma_A(Z)=<JV,R(\pi_1(A),Z)J>,\quad Z\in T_pM.$$
Then
$$
[A^h,{\cal I}\widetilde W]_K=2\gamma_A^h.
$$

Computations in local coordinates  involving (\ref{cal I/hor}), (\ref{cal I/ver}),
(\ref{eq Xh-y-z}) and (\ref{eq der a-b})  show that
$$
 [{\cal I}A^h,\widetilde W]_K=-2(U(A))^{h}_{K}.
$$
and
$$
[{\cal I}A^h,{\cal I}\widetilde W]_K=-2((IU)(A))^{h}_{K}+ 2\gamma_{IA}^h.
$$

It follows that
$$N^{{\cal I}}(A^h,W)=\frac{1}{2}N^{{\cal I}}(A^h,\widetilde W)_K= -{\cal
J}\gamma_A^h+\gamma_{IA}^h.$$

  Let $\{e_1,e_2\}$ be a bases of $T_pM$ and denote by $\{\eta_1,\eta_2\}$ its dual
bases. Then $Q_1=e_1+\eta_1$, $Q_2=e_2+\eta_2$, $Q_3=e_1-\eta_1$, $Q_4=e_2-\eta_2$
constitute an orthonormal bases of $T_pM\oplus T_pM^{\ast}$ yielding its canonical
orientation. According to Example 5, every generalized complex structure $J\in
G^{-}(T_pM)$ is given by
$$
\begin{array}{lll}
Q_1\to y_1Q_2+y_2Q_3+y_3Q_4,& & Q_2\to -y_1Q_1+y_2Q_4-y_3Q_3 \\
Q_3\to -y_1Q_4+y_2Q_1-y_3Q_2,& & Q_4\to y_1Q_3+y_2Q_2+y_3Q_1,
\end{array}
$$
where $y_1^2-y_2^2-y_3^2=1$, $y_1,y_2,y_3\in {\Bbb R}$. Then
$$
\begin{array}{c}
{\cal J}\gamma_A^h=\gamma_A(e_1)(J\eta_1)^h+\gamma_A(e_2)(J\eta_2)^h=\\[4pt]
-(y_1+y_3)\gamma_A(e_2)e_1^h+(y_1+y_3)\gamma_A(e_1)e_2^h-
y_2\gamma_A(e_1)\eta_1^h-y_2\gamma_A(e_2)\eta_2^h.
\end{array}
$$
Therefore the identity $N^{{\cal I}}(A^h,W)=0$ implies $\gamma_A(e_1)=\gamma_A(e_2)=0$,
i.e. $\gamma_A=0$. This proves statement $(i)$. The proof of $(ii)$ is similar.
\end{proof}

\smallskip

Now suppose that $R(X,Y)I=0$ for every generalized complex structure $I\in
G^{+}(T_pM)$, $X,Y\in T_pM$ being fixed. Take a basis $\{e_1,e_2\}$ of $T_pM$, denote
by  $\{\eta_1,\eta_2\}$ its dual bases and set $Q_1=e_1+\eta_1$, $Q_2=e_2+\eta_2$,
$Q_3=e_1-\eta_1$, $Q_4=e_2-\eta_2$. Then every $I$ is given by (see Example 5)
$$
\begin{array}{lll}
Q_1\to x_1Q_2+x_2Q_3+x_3Q_4,& & Q_2\to -x_1Q_1-x_2Q_4+x_3Q_3 \\
Q_3\to -x_1Q_4+x_2Q_1+x_3Q_2,& & Q_4\to -x_1Q_3-x_2Q_2+x_3Q_1,
\end{array}
$$
where $x_1^2-x_2^2-x_3^2=1$, $x_1,x_2,x_3\in {\Bbb R}$. The identity $R(X,Y)I=0$
implies $<R(X,Y)Ie_1,\eta_k>+<R(X,Y)e_1,I\eta_k>=0$, $k=1,2$, which is equivalent to
$$
(x_1+x_3)\eta_1(R(X,Y)e_2)+(x_1-x_3)\eta_2(R(X,Y)e_1)=0,
$$
$$
2x_2\eta_2(R(X,Y)e_2)-(x_1+x_3)\eta_1(R(X,Y)e_1)+(x_1+x_3)\eta_2(R(X,Y)e_2)=0.
$$
It follows that $R(X,Y)I=0$ for every $I$ if and only if $R(X,Y)=0$.

  It is also easy to see that $R(X,Y)J=0$ for every $J\in G^{+}(T_pM)$
if and only if $\eta_1(R(X,Y)e_1)+\eta_2(R(X,Y)e_2)=0$.

  Thus if the structures ${\cal I}$ and ${\cal J}$ are both integrable, then the
connection $\nabla$ is flat. The converse is also true as the following result shows.

\begin{tw}\label {Integrability}
Let $M$ be a $2$-dimensional manifold and $\nabla$ a torsion-free  connection on $M$.
Then the generalized almost complex structures ${\cal I}$ and ${\cal J}$  induced by
$\nabla$ on the twistor space ${\cal P}$ are both integrable if and only if the
connection $\nabla$ is flat.
\end{tw}

\begin{proof}
 Since the structures ${\cal I}$ and ${\cal J}$ on ${\cal V}\oplus{\cal V}^{\ast}$
are induced by complex structures on the fibres of ${\cal P}$ the Nijenhuis tensors of
${\cal I}$ and ${\cal J}$ vanish on ${\cal V}\oplus{\cal V}^{\ast}$. Thus, in view of
Propositions~\ref{Nijenhuis-hor} and \ref{Nijenhuis-hor-ver}, we have to consider these
tensors only on ${\cal H}\times{\cal V}^{\ast}$.

  Suppose that the connection $\nabla$ is flat. Let $K=(I,J)\in{\cal P}$. Fix bases $\{U_1,U_2={\cal K}^{+}U_1\}$ of ${\cal V}_{I}{\cal
G}^{+}$ and $\{V_1,V_2={\cal K}^{-}V_1\}$ of ${\cal V}_{J}{\cal G}^{-}$. Take sections
$a_1$ and $b_1$ of $A(M)$ near the point $p=\pi(K)$ such that $a_1(p)=U_1$,
$b_1(p)=V_1$ and $\nabla a_1|_p=\nabla b_1|_p=0$. Define vertical vector fields
$\widetilde{a_1}$ and $\widetilde{b_1}$ on ${\cal G}^{+}$ and ${\cal G}^{-}$ by means
of (\ref{eq def-ver-ab}). Set $\widetilde{a_2}={\cal K}^{+}\widetilde{a_1}$,
$\widetilde{b_2}={\cal K}^{-}\widetilde{b_1}$. Then
$\{\widetilde{a_1},\widetilde{a_2}\}$ and $\{\widetilde{b_1},\widetilde{b_2}\}$ are
frames of the vertical bundles ${\cal V}{\cal G}^{+}$ and ${\cal V}{\cal G}^{-}$ near
the points $I$ and $J$, respectively. Denote by $\{\alpha_1,\alpha_2\}$ and
$\{\beta_1,\beta_2\}$ the dual frames of $\{\widetilde{a_1},\widetilde{a_2}\}$ and
$\{\widetilde{b_1},\widetilde{b_2}\}$. Set $\widetilde{W_i}=(\widetilde{a_i},0)$,
$\gamma_i=(\alpha_i,0)$ and $\widetilde{W_{i+2}}=(0,\widetilde{b_i})$,
$\gamma_{i+2}=(0,\beta_i)$ for $i=1, 2$. Then $\{\widetilde {W_1},\widetilde
{W_2},\widetilde {W_3},\widetilde {W_4}\}$ is a frame of the vertical bundle ${\cal V}$
of ${\cal P}$ near the point $K$ and $\{\gamma_1,\gamma_2,\gamma_3,\gamma_4\}$ is its
dual frame. We have $\gamma_2={\cal I}\gamma_1$, ${\cal I}\gamma_3=\beta_2^{\sharp}$,
${\cal I}\gamma_4=-\beta_1^{\sharp}$. If $A\in T_pM\oplus T_p^{\ast}M$, then ${\cal
I}N^{\cal I}(A^h,\gamma_3)=-N^{\cal I}(A^h,{\cal I}\gamma_3)=-N^{\cal
I}(A^h,\beta_2^{\sharp})=0$ by Proposition~\ref{Nijenhuis-hor-ver}. Hence $N^{\cal
I}(A^h,\gamma_3)=0$. Similarly, $N^{\cal I}(A^h,\gamma_4)=0$.

As in the proof of Proposition~\ref{Nijenhuis-hor-ver}, it is not hard to see that
$[\pi_1(A^h),\widetilde{W_r}]_K=0$,  $r=1,...,4$, $[\pi_1({\cal
I}A^h),\widetilde{W_i}]_K=-(\pi_1(I\widetilde{ a_i}(A)))^h_K$ and $[\pi_1({\cal
I}A^h),\widetilde{W_{i+2}}]_K=0$, $i=1, 2$. In particular
$[\pi_1(A^h),\widetilde{W_r}]_K$   and $[\pi_1({\cal I}A^h),\widetilde{W_r}]_K$ are
horizontal vectors for every $r=1,...,4$. It follows, in view of (\ref{eq 3.2}) and
Lemma~\ref{brackets}$(i)$, that for every $Z\in T_pM$, $r=1,...,4$ and $s=1, 2$
$$
\begin{array}{c}
({\cal L}_{\pi_1(A^h)}\gamma_s)_K(Z^h+\widetilde W_r)=-\alpha_s(R(\pi_1(A),Z)I)=0, \\[4pt]
({\cal L}_{\pi_1({\cal I}A^h)}\gamma_s)_K(Z^h+\widetilde
W_r)=-\alpha_s(R(\pi_1(IA),Z)I)=0
\end{array}
$$
since the connection $\nabla$ is flat. This implies $N^{\cal I}(A^h,\gamma_s)_K=0$ for
$s=1, 2$.

It follows that $N^{\cal I}(A^h,\Theta)_K=0$  for every $\Theta\in{\cal V}^{\ast}_K$.
Similarly, $N^{\cal J}(A^h,\Theta)_K=0$.
\end{proof}

\smallskip

 Denote by $(g,J_{+},J_{-},b)$ the data on ${\cal P}$ determined by the almost
generalized K\"ahler structure  $\{{\cal I}^{\nabla},{\cal J}^{\nabla}\}$ as described
in \cite{Gu}. It is not hard to see that the metric $g$, the almost complex structures
$J_{\pm}$ and the $2$-form $b$ are given as follows.  Let $K=(I,J)\in{\cal P}$, $X,Y\in
T_{\pi(K)}M$, $W=(U,V)\in{\cal V}_K$. Let $\{e_1,e_2\}$ be a local frame of $TM$ near
the point $\pi(K)$ and denote by $\{\eta_1,\eta_2\}$ its dual co-frame. Define
endomorphisms $I_r$, $J_s$, $r,s=1,2,3$, by means of $e_1,e_2,\eta_1,\eta_2$ as in
Example 5. Then $I=\sum_{r}x_rI_r$, $J=\sum_s y_sJ_s$ with $x_1^2-x_2^2-x_3^2=1$,
$y_1^2-y_2^2-y_3^2=1$.  Let $X=X_1e_1+X_2e_2$, $Y=Y_1e_1+Y_2e_2$. Then
$$
\begin{array}{c}
g(X^h,Y^h)_K=\\
\displaystyle{\frac{1}{y_1+y_3}}[(x_1+x_3)X_1Y_1-x_2(X_1Y_2+X_2Y_1)+(x_1-x_3)X_2Y_2],\\[6pt]
g(X^h,W)_K=0,\quad g|({\cal V}_K\times {\cal V}_K)=h.
\end{array}
$$
$$
\begin{array}{lll}
J_{+}X^h_K=(IX)^h_K, & &  J_{-}X^h_K=(IX)^h_K,  \\
J_{+}(U,V)=(I\circ U,J\circ V), & & J_{-}(U,V)=(I\circ U,-J\circ V).
\end{array}
$$
$$
\begin{array}{c}
b(X^h,Y^h)_K=\displaystyle{\frac{y_2}{y_1+y_3}}(X_1Y_2-X_2Y_1),\\[6pt]
b(X^h,W)_K=0,\quad b|({\cal V}_K\times {\cal V}_K)=0.
\end{array}
$$
In particular, the almost complex structures $J_{+}$ and $J_{-}$ commutes and
$J_{+}\neq\pm J_{-}$.

Computations similar to that above show that the almost complex structures $J_{\pm}$
are both integrable for any torsion-free connection $\nabla$. Denote by $\omega_{\pm}$
the K\"ahler form of the Hermitian structure $(g,J_{\pm})$ on ${\cal P}$. Then
$$
\begin{array}{c}
\omega_{\pm}(X^h,Y^h)_K= (y_1+y_3)^{-1}(X_1Y_2-X_2Y_1),\quad \omega_{\pm}(X^h,W)_K=0,
\\[6pt]
 \omega_{\pm}(W,W^{'})=h(I\circ U,U^{'})\pm h(J\circ V,V^{'}),
\mbox { where } W^{'}=(U^{'},V^{'})\in{\cal V}_K.
\end{array}
$$
Set $V=\sum_s v_sJ_s$. Then we easily obtain that
$$
\begin{array}{lll}
3d\omega_{\pm}(X^h,Y^h,W)_K&=&-(v_1+v_3)(y_1+y_3)^{-2}((X_1Y_2-X_2Y_1)\\[6pt]
&+&h(R(X,Y)I,I\circ U)\pm h(R(X,Y)J,J\circ V)
\end{array}
$$
in view of (\ref{eq 3.2}) and the fact that $[X^h,W]_K$ and $[Y^h,W]_K$ are vertical
vectors. Moreover
$$
\begin{array}{lll}
h(R(X,Y)J,J\circ V)=-<R(X,Y)J,J\circ V>=\\[4pt]
2(y_1+y_3)[y_2(v_1-v_3)+v_2(y_1-y_3)][\eta_1(R(X,Y)e_1)+\eta_2(R(X,Y)e_2)]
\end{array}
$$
Thus putting $y_1=2$, $y_2=0$, $y_3=\sqrt 3$, $U=0$, $v_1=\sqrt 3$, $v_2=0$, $v_3=2$ we
see that $d\omega_{\pm}(X^h,Y^h,W)\neq 0$. Therefore the structure $(g,J_{\pm})$ is not
K\"ahlerian.

\bigskip

\noindent {\bf Acknowledgement}: We would like to thank W. Goldman for sending us a
proof of the classification of complete affine $2$-dimensional manifolds.

\end{document}